\documentclass[11pt,intlimits]{amsart}
\begin{document}
\newtheorem{theorem}{Theorem}[section]
\newtheorem{proposition}[theorem]{Proposition}
\newtheorem{lemma}[theorem]{Lemma}
\newtheorem{corollary}[theorem]{Corollary}
\newtheorem{definition}[theorem]{Definition}
\newtheorem{conjecture}[theorem]{Conjecture}

\numberwithin{equation}{section}
\renewcommand{\theequation}{\thesection.\arabic{equation}}

\newcommand{\princEigen}{\mathfrak{E}}
\newcommand{\bp}{\mathfrak{b}}
\newcommand{\xib}{\boldsymbol{\xi}}
\newcommand{\me}{\mathrm{e}} 
\newcommand{\dif}{\mathrm{d}}
\newcommand{\rdm}{\mathfrak{m}}
\newcommand{\jactr}{\nabla\! {\mathbf X}^{\dagger}}
\newcommand{\jac}{\nabla\! {\mathbf X}}
\newcommand{\dersigmatr}{\nabla\!\mbox{\boldmath$\sigma$}_j^{\dagger}}
\newcommand{\dersigma}{\nabla\!\mbox{\boldmath$\sigma$}_j}
\newcommand{\derdrift}{\nabla\!{\mathbf m}}
\newcommand{\derdriftstrat}{\nabla\!{\mathbf q}}
\newcommand{\forfilt}[2]{\ensuremath{\mathbb{B}_{#1}^{#2}}}  
\newcommand{\backfilt}[2]{\ensuremath{\mathbb{B}_{#1}^{#2}}} 
\newcommand{\filt}[1]{\forfilt{}{#1}}
\newcommand{\bdown}{\ensuremath{b^{\downarrow}_T}}            
\newcommand{\bback}{B}
\newcommand{\lifetimesharp}{\ensuremath{\mathfrak{e}}}        
\newcommand{\invmeas}{\ensuremath{\dif\Pi}}                   
\newcommand{\xinf}{\ensuremath{\mathfrak{X}}}                 
\newcommand{\real}[1]{\ensuremath{\mathbb{R}^{#1}}}           

\newcommand{\gup}{\ensuremath{\mathcal{G}}}                      
\newcommand{\gsharp}{\ensuremath{\mathcal{G}^\sharp}}        
\newcommand{\gupAdj}{\ensuremath{\mathfrak{L}}}                      
\newcommand{\gsharpAdj}{\ensuremath{\mathfrak{L}^\sharp}}        

\newcommand{\subHarmonicU}{\ensuremath{u}}           

\newcommand{\trace}{\mathtt{tr}}                              
\newcommand{\strat}{\!\bullet\!}
\newcommand{\Czero}{\ensuremath{\mathcal{C}}}       
\newcommand{\Czerob}{\ensuremath{\mathcal{C}_b}}    
\newcommand{\Ctwo}{\ensuremath{\mathcal{C}^2}}    
\newcommand{\Ctwob}{\ensuremath{\mathcal{C}^2_b}} 
\newcommand{\smooth}{\ensuremath{\mathcal{C}^\infty}} 
\newcommand{\Cvan}{\ensuremath{\mathcal{C}_\infty}} 

\newcommand{\proup}[1]{\ensuremath{X_{#1}}}                
\newcommand{\proupb}[1]{\ensuremath{\mathbf{X}_{#1}}}      
\newcommand{\prodown}[1]{\ensuremath{X^{\downarrow}_{#1}}} 
\newcommand{\prodownb}[1]{\ensuremath{\mathbf{X}^{\downarrow}_{#1}}} 
\newcommand{\prointer}{\ensuremath{X^T}}                  
\newcommand{\prointermediate}[1]{\ensuremath{X^{T}_{#1}}}    
\newcommand{\prosharp}[1]{\ensuremath{X^{\sharp}_{#1}}}        
\newcommand{\prosharpb}[1]{\ensuremath{\mathbf{X}^{\sharp}_{#1}}}   
\newcommand{\prosharpinv}[1]{\ensuremath{X^{\sharp -1}_{#1}}}   
\newcommand{\proback}[1]{\ensuremath{\mathfrak{z}_{#1}}}             

\newcommand{\flowup}[1]{\ensuremath{X_{#1}}}                   
\newcommand{\flowdown}[1]{\ensuremath{X^{\downarrow}_{#1}}}    
\newcommand{\flowintermediate}[1]{\ensuremath{X^{T}_{#1}}}     
\newcommand{\flowsharp}[1]{\ensuremath{X^{\sharp}_{#1}}}       
\newcommand{\flowsharpinv}[1]{\ensuremath{X^{\sharp -1}_{#1}}} 
\newcommand{\invmap}[1]{\ensuremath{\phi_{#1}}}

\bibliographystyle{plain}

\title[On the Lyapunov exponent]{On the Lyapunov exponent of a multidimensional stochastic flow}
\author[M. Baldini]{Michele Baldini$^\dag$}
\thanks{$^\dag$Merrill Lynch, Global Equity Linked Products, New York. This paper was written while the author was visiting Northwestern University and it represents only his personal opinions and not those of Merrill Lynch, its subsidiaries or affiliates.}

\begin{abstract}
Let $\proup{t}$ be a reversible and positive recurrent diffusion in \real{d} described by
\begin{equation}\nonumber
X_t=x+\sigma\, b(t)+\int_0^tm(X_s)\dif s,
\end{equation}
where the diffusion coefficient $\sigma$ is a positive-definite matrix and the drift $m$ is a smooth function.
Let $\proup{t}(A)$ denote the image of a compact set $A\subset\real{d}$ under the stochastic flow generated by $\proup{t}$.
If the divergence of the drift is strictly negative, there exists a set of functions $u$ 
such that 
\[\lim_{t\to\infty}
\int_{\proup{t}(A)}\subHarmonicU(x)\dif x=0\quad\mbox{a.s.}
\]
A characterization of the functions $u$ is provided, as well as lower and upper bounds for the exponential rate of convergence.
\end{abstract}
\maketitle

\section{Introduction}
In \cite{hasminskii_focusing} Has'minskii proved that any solution of a one-dimensional stochastic differential equation
describing a positive recurrent diffusion \proup{t} is stable with probability one in the metric given by the natural scales.
More precisely, let $\proup{t}$ be a diffusion with smooth coefficients
\begin{equation}
X_t=x+\int_0^t\sigma(X_s)\dif b(s)+\int_0^tm(X_s)\dif s.
\end{equation}
Then, for every pair of initial conditions $y,z\in\real{}$
\begin{equation}
{\bf P}\left[\lim_{t\to\infty}\big|s(\proup{t}(y))-s(\proup{t}(z))
  \big|=0   \right]=1,
\end{equation}
where $s$ is the scale function. A different proof of this stability is offered in \cite{Baldini_thesis}, together with the result that the pathwise rate of convergence is exponentially fast, i.e.
\begin{equation}
{\bf P}\left[\lim_{t\to\infty}\frac{1}{t}\ln\big|s(\proup{t}(y))-s(\proup{t}(z))\big|=-2\int_\real{}\frac{m^2}{\sigma^2}\dif\Pi\right]=1.
\end{equation} 
Here $\dif\Pi$ is the corresponding invariant measure.

This paper proceeds in the same direction and tackles the multidimensional problem. As expected, technical difficulties make the task quite challenging and more restrictive conditions on the coefficients are introduced. The main result is: Let $\proup{t}$ be the stochastic flow generated by a reversible and positive recurrent diffusion in $\real{d}$. If the diffusion coefficient is a positive-definite matrix and the drift is a $C^1$ function subject to certain boundness conditions on the derivatives, then there exists a non-empty set of smooth functions $\subHarmonicU$ such that, for every compact set $A\subset\real{d}$, the $\subHarmonicU$-volume of the flow $\proup{t}(A)$ tends almost surely to zero as $t\uparrow\infty$. In addition, the upper and lower bounds for the corresponding Lyapunov exponent are provided. 

The paper is organized as follows: section 2 introduces the notation and the necessary background. Section 3 paves the road with some introductory results. The following two sections 
are devoted to the main theorems:  section 4 studies the asymptotic behavior of the flow and section 5 discusses the Lyapunov exponent. An example is illustrated in section 6 and section 7 concludes. 

\section{Notation and preliminaries}
Throughout this article it is assumed that:
$d\ge 2$ is the dimension of the phase space $\real{d}$, $\mathcal{B}$ is the $\sigma$-algebra of Borel sets of $\real{d}$ with the implicit assumption that 
every subset in $\real{d}$ belongs to $\mathcal{B}$.
Given $D, D'\subseteq\real{d}$, $D'\subset\subset D$ denotes a bounded subset such that $\bar{D}'\subset D$. $A$ is always a compact subset and $\mathrm{Leb}(A)$ is the Lebesgue measure of $A$. 

Let $D\subseteq\real{d}$. A function $f$ in $D$ is called \textit{uniformly H\"{o}lder continuous with exponent $\alpha\in(0,1]$} if
\begin{equation}
\sup_{\substack{x,y\in D\\ x\neq y}}\frac{|f(x)-f(y)|}{|x-y|^\alpha}<\infty.
\end{equation}
If the same inequality holds for every $D'\subset\subset D$, then $f$ is called \textit{locally H\"{o}lder continuous}. A function $f$ is called \textit{uniformly (locally) Lipschitz continuous} if it is uniformly (locally) H\"{o}lder continuous with exponent 1. Define the H\"{o}lder space $C^{k,\alpha}(\bar{D})(C^{k,\alpha}(D))$ as the space of functions whose $k$th order partial derivatives are uniformly (locally) H\"{o}lder continuous with exponent $\alpha\in(0,1]$ in $D$.

Let $(\Omega, \mathcal{F}, {\bf P})$ be the canonical Wiener space of the $d$-dimensional Brownian motion with paths $b(t),\,t\ge 0$ and let $\theta_t$ be the measure-preserving shift operator on the space of Brownian paths, i.e. $\theta_t b(s):=b(s+t)$. Let $\sigma$ be a $d\times d$ positive-definite matrix and $a:=\sigma\sigma^t$. Let $m:\real{d}\to\real{d}$ be a drift coefficient such that:
\begin{enumerate}
\renewcommand{\theenumi}{\roman{enumi}}
\item $m\in C^{1,1}(\real{d})$.
\item There exists a potential $\psi:\real{d}\to(0,\infty)$ that satisfies
\begin{eqnarray}
& & m=-a\nabla\ln{\psi}\label{reversibility},\\[2mm]
& & \int_{\real{d}}\frac 1{\psi^2}=1.\label{finite_invariant_measure}
\end{eqnarray}
\item There exists a constant $\lambda_c>0$ such that
\begin{equation}
(\nabla\cdot m)(x)\le -\lambda_c\quad\mbox{for all }x\in\real{d}.\label{negativeDivergence}
\end{equation}
\end{enumerate}
Roughly speaking, condition (iii) says that the drift resembles a restoring force \textit{\`{a} la 
Ornstein-Uhlenbeck}, i.e. $m(x)=-x$.
The reader will also notice that the potential is defined modulo a positive multiplicative constant, so to assure \eqref{finite_invariant_measure} it suffices to verify that $\psi^{-2}$ is integrable.

Let $\proup{t}(x,b)$ be the flow determined by the It\^{o} stochastic differential equation (SDE)
\begin{equation}\label{recurrent_diffusion}
X_t=x+\sigma\, b(t)+\int_0^tm(X_s)\dif s.
\end{equation}
Here ``flow" means that $\proup{t}(x,\cdot)$ is a function of a single Brownian motion, the same for every $x$, so that for every $t\ge 0$, $\proup{t}(x,\cdot)$ is a 1-1 map from $\real{d}$ into $\real{d}$. More precisely, after a null set in $\Omega$ has been weeded out, the family of maps $\{\proup{t}(x,b):t\ge 0\}$ defines a stochastic flow of global diffeomorphisms, i.e. for every Brownian path $b$
\begin{enumerate}
\renewcommand{\theenumi}{\roman{enumi}}
\item $\real{d}\ni x\mapsto\proup{t}(x,b)$ is a diffeomorphism in $\real{d}$ for all $t\ge0$ and the inverse map is continuous in $t$ and smooth in $x$.
\item $\proup{t+\tau}(x,b)=\proup{t}(\cdot,\theta_{\tau}b)\circ\proup{\tau}(x,b)$ for all non-negative $t,\tau$ and $x\in\real{d}$.
\end{enumerate}
For an exceptional presentation of this topic and for the proof that uniformly Lipschitz continuous coefficients generate a flow of global diffeomorphisms, the reader should consult \cite{kunita}, section 4.7.

The infinitesimal generator $\gup$ of the diffusion introduced above can now be worked out. For any $f\in C^2(\real{d})$, 
\begin{equation}
\gup f = \frac 1{2}\psi^2\nabla\cdot(\psi^{-2}a\nabla f) = \frac 1{2}\nabla\cdot a\nabla f +m\cdot\nabla f,\label{g_up}
\end{equation}
whereas the formal adjoint operator is 
\begin{equation}
\gupAdj f = \frac 1{2}\nabla\cdot[\psi^{-2}a\nabla(\psi^{2}f)]=\frac 1{2}\nabla\cdot a\nabla f -m\cdot\nabla f-(\nabla\cdot m)f.\label{g_up_adj}
\end{equation}

Thanks to \eqref{g_up}-\eqref{g_up_adj}, 
the meaning of \eqref{reversibility} and \eqref{finite_invariant_measure} becomes clear: the first guarantees that $\gup$ is self-adjoint in $L^2(\frac{\dif x}{\psi^2(x)})$ (this is the {\it reversibility} condition in \cite{Durrett_reversible_diffusion}), the latter says that 
the invariant density measure $\frac 1{\psi^2}$ of \proup{t} is indeed an invariant probability density. 
For general questions on the existence and uniqueness of solutions of SDE, we refer the reader to \cite{mckean_stoch_int}.

An allied flow is now presented: Let $\prosharp{t}(x,b)$ be the flow generated by the 
SDE 
\begin{equation}\label{sharp_transient_diffusion}
X^{\sharp}_t=x+\sigma\, b(t)-\int_0^tm(X^{\sharp}_s)\dif s
\end{equation}
with reversed drift $-m$ in place of $m$. Since the regularity conditions on the coefficients are the same, the family of maps $\{\prosharp{t}(x,b):t\ge 0\}$ defines a stochastic flow of global diffeomorphisms of $\real{d}$. However, it will be shown that the ergodic property of the corresponding diffusions $\prosharp{t}$ is diametrically opposite to the one displayed by $\proup{t}$. 

The diffusion $\prosharp{t}$ is not the main subject of this paper, however it is intimately related to $\proup{t}$ (for more of this see \cite{Baldini_thesis}).
In the next section the ergodicities of both processes are studied, in conjunction with the harmonic property of $\prosharp{t}$'s infinitesimal generator: it turns out that this operator enjoys a well-defined criticality condition which will be the key for the proof of the main theorems in the conclusive sections.

For convenience of the reader, the generator and the adjoint operator are spelled out in full form: 
\begin{eqnarray}
& & \gsharp f = \frac 1{2}\psi^{-2}\nabla\cdot(\psi^2 a\nabla f) = \frac 1{2}\nabla\cdot a\nabla f-m\cdot\nabla f,\label{g_sharp}\\[2mm]
& & \gsharpAdj f = \frac 1{2}\nabla\cdot[\psi^2 a \nabla(\psi^{-2} f)]=\frac 1{2}\nabla\cdot a\nabla f +m\cdot\nabla f+(\nabla\cdot m)f\label{g_sharp_adj}.
\end{eqnarray}

\section{Ergodic properties of \proup{t} and \prosharp{t}}
A criterion of R. Pinsky can be applied to study recurrence and transience in terms of integral tests.

\begin{theorem}[\cite{pinsky}, pp. 258]\label{pinsky_criteria}
Let $L= \frac 1{2}\nabla\cdot a\nabla+a\nabla Q\cdot\nabla$ on $\real{d},\ d\ge 2$ satisfying the following conditions on every $D\subset\subset\real{d}$: $a_{ij}(x)\in C^{1,\alpha}(\bar{D})$, $Q(x)\in C^{2,\alpha}(\bar{D})$ and $\sum_{i,j=1}^d a_{ij}(x)v_iv_j>0\mbox{, for all }v\in\real{d}-\{0\}\mbox{ and all }x\in\bar{D}$. Define
\begin{eqnarray}
E_1(x) & = & \frac{(a(x)x,x)}{|x|^2}e^{2Q(x)},\\[2mm]
E_2(x) & = & \frac{|x|^2}{(a^{-1}(x)x,x)}e^{2Q(x)},\\[2mm]
\hat{E}_1(r) & = & \int_{S^{d-1}}E_1(r\theta)\dif \theta,
\end{eqnarray}
where $\dif\theta$ is Lebesgue measure on $S^{d-1}$.
\begin{enumerate}\renewcommand{\theenumi}{\roman{enumi}}
\item If 
\begin{equation}
\int_1^\infty r^{1-d} \frac 1{\hat{E}_1(r)}\dif r = \infty,
\end{equation}
then the diffusion corresponding to $L$ is recurrent.
\item If 
\begin{equation}
\int_1^\infty r^{1-d} \frac 1{E_2(r\theta)}\dif r < \infty,
\end{equation}
on a subset of $S^{d-1}$ with positive Lebesgue measure, then the diffusion corresponding to $L$ is transient.
\end{enumerate}
\end{theorem}

\begin{corollary}
\proup{t} is positive recurrent and \prosharp{t} is transient.
\end{corollary}
\begin{proof} Let $\lambda_{m}$ and $\lambda_M$ be respectively the smallest and largest eigenvalue of $a$. Then $0<\lambda_m\le \lambda_M <\infty$. A straightforward computation gives $\lambda_m |x|^2 \le (ax,x) \le \lambda_M |x|^2$ and 
$\frac {|x|^2}{\lambda_M} \le (a^{-1}x,x) \le \frac {|x^2|}{\lambda_m}$. Now, by virtue of the integrability of $\psi^{-2}$, it follows that 
\begin{eqnarray}
& & \int_1^\infty \frac{r^{1-d}}{\int_{S^{d-1}}\psi^{-2}(r\theta)\dif\theta}\dif r = \infty\quad\mbox{and}\\[4mm]
& & \int_1^\infty \frac{r^{1-d}}{\psi^{2}(r\theta)}\dif r<\infty\mbox{ for almost all } \theta\in S^{d-1}.
\end{eqnarray}
A direct application of theorem \ref{pinsky_criteria} establishes recurrence and transience for $\proup{t}$ and $\prosharp{t}$, respectively. What remains is to verify the \textit{positiveness} of \proup{t}. This follows from the observation that $\psi^{-2}$ is the invariant probability density of $\proup{t}$ (see discussion in \cite{pinsky}, chapter 4.9).
\end{proof}

Some spectral properties of the operator $\gsharpAdj$ are explored next.
By the transience property, both $\gsharp$ and $\gsharpAdj$ possess a Green function (\cite{pinsky}, Th. 4.3.3). This implies that for each non-negative $f\in C_0^1(\real{d})$ which is non identically zero, there exist positive solutions $u\in C^{2,1}(\real{d})$ of $\gsharpAdj u=-f$. The smallest such solution is given by $u_0(x)=\int_{\real{d}}G(x,y)f(y)\dif y$ where $G$ is the Green function of $\gsharpAdj$ (\cite{pinsky}, Th. 4.3.8). Define the cone of all superharmonic functions as follows:
\begin{equation} 
\mathcal{S}=\{0 < u\in C^{2,1}(\real{d}):0\ge \gsharpAdj u\in  C_0^1(\real{d}),\,\gsharpAdj u\mbox{ non identically zero}\}.
\end{equation}
Denote $\princEigen$ the principal eigenvalue of $\gsharpAdj$, i.e. the supremum of the spectrum of the corresponding self-adjoint extension. Using Theorem 4.4.5 in \cite{pinsky} and the transience property it follows that 
\begin{equation}
\princEigen=\inf_{\substack{u\in C^{2,1}(\real{d})\\ u>0}}\ \sup_{x\in\real{d}}\frac{\gsharpAdj u}{u}(x)\leq 0.
\end{equation}
Thanks to the negative divergence of the drift, one can prove more.
\begin{theorem}\label{princEigen_is_negative}
The principal eigenvalue $\princEigen$ is strictly negative.
\end{theorem}
\begin{proof}
As stated above, for each non-negative $f\in C_0^1(\real{d})$ which is not identically zero there always exists a positive $u\in C^{2,1}(\real{d})$ such that $\gsharpAdj u =-f\le 0$.
Now, with $g\in C^2(\real{d})$, the following identity holds:
\begin{equation}\label{operatorOnSquareFunction}
\gsharpAdj g^2 = 2 g\, \gsharpAdj g -g^2 \nabla\cdot m +\nabla g\cdot a\nabla g.
\end{equation}
Since the square root does not spoil the regularity of a positive function, apply 
\eqref{operatorOnSquareFunction} to $g=\sqrt{u}$ obtaining
$2\sqrt{u}\, \gsharpAdj \sqrt{u} =  \gsharpAdj u +u \nabla\cdot m -\nabla \sqrt{u}\cdot a\nabla \sqrt{u}$. Using \eqref{negativeDivergence} and the positive definiteness of $a$, we can write
\begin{equation}
\frac{\gsharpAdj \sqrt{\subHarmonicU}}{\sqrt{\subHarmonicU}}\le -\frac{\lambda_c}{2}.
\end{equation}
The proof is completed observing that 
\begin{equation}
\princEigen\le\sup_{x\in\real{d}}\frac{\gsharpAdj \sqrt{\subHarmonicU}}{\sqrt{\subHarmonicU}}(x)\le -\frac{\lambda_c}{2}.
\end{equation}
\end{proof}

\section{Main results}\label{focusingflow}
Denote $\proup{t}^i$ the $i$-th component of \proup{t} and $m^i_j$ the $j$-th derivative of the $i$-th component of $m$. Denote $|D\proup{t}|:=|D\proup{t}(x,b)|$ the Jacobian of the map $x\mapsto\proup{t}(x,b)$. 
\begin{lemma}\label{lemma_uXdet(jac)}
Let $0<u\in C^2(\real{d})$. Then
\begin{equation}\nonumber
u(\proup{t})|D\proup{t}|=u(x)\mathrm{exp}\left(\int_0^t\frac{\nabla u(\proup{\tau})}{u(\proup{\tau})}\cdot\sigma\dif b -\frac{1}{2}\int_0^t \frac{(\nabla u\cdot a\nabla u)(\proup{\tau})}{u^2(\proup{\tau})}\dif \tau+\int_0^t\frac{(\gsharpAdj u)(\proup{\tau})}{u(\proup{\tau})}\dif \tau  \right).
\end{equation}
\end{lemma}
\begin{proof}
From standard results \cite{kunita}, the gradient flow satisfies $\frac{\partial\proup{t}^i}{\partial x^j}=\delta^{ij}+\sum_{k=1}^d\int_0^t m^i_k(\proup{\tau})\frac{\partial\proup{\tau}^k}{\partial x^j}\dif\tau$. Simple differentiation
gives
\begin{eqnarray}
\dif |D\proup{t}| & = & \sum_{\mathcal{P}}(-1)^p\left(\sum_{l=1}^d \frac{\partial\proup{t}^1}{\partial x^{i_1}}\dots \sum_{k=1}^d m^l_k(\proup{t})\frac{\partial\proup{t}^k}{\partial x^{i_l}}\dots \frac{\partial\proup{t}^d}{\partial x^{i_d}}    \right)\dif t\nonumber \\[3mm]
& = & \sum_{k,l=1}^d m^l_k(\proup{t})\left( \sum_{\mathcal{P}}(-1)^p \frac{\partial\proup{t}^1}{\partial x^{i_1}} \dots \frac{\partial\proup{t}^k}{\partial x^{i_l}}\dots \frac{\partial\proup{t}^d}{\partial x^{i_d}}      \right)\dif t \\[3mm]
& = & \nabla\cdot m(\proup{t})|D\proup{t}|\dif t\nonumber,
\end{eqnarray}
where $\mathcal{P}$ is the set of all permutations $(1,\dots,d)\to(i_1,\dots,i_d)$ and $p$ is the sign of each permutation. Now use It\^{o}'s lemma on $u(\proup{t})|D\proup{t}|$ to show  
\begin{multline}\label{M_supermartingale}
\dif \Big[u(\proup{t})|D\proup{t}|\Big]\\
=|D\proup{t}|\left[\nabla u(\proup{t})\cdot \sigma\dif b + \left(\frac{1}{2}\nabla\cdot a\nabla u + m\cdot\nabla u + (\nabla\cdot m)u\right)(\proup{t})\dif t  \right].
\end{multline}
The proof is completed observing that the operator acting on $u$ on the RHS of 
\eqref{M_supermartingale} is the formal adjoint of $\gsharp$.
\end{proof}

Next we turn to the analysis of the volume of $\proup{t}(A)$, i.e. the image of a compact set $A\subset\real{d}$ under the flow $\proup{t}(\cdot,b)$.

\begin{theorem}\label{volume_flow_is_supermartingale}
Let $0<\subHarmonicU\in C^2(\real{d})$ such that $\gsharpAdj\subHarmonicU\le 0$. Then, 
\begin{equation}
\int_{\proup{t}(A)}\!\!\!\subHarmonicU:=\int_{\proup{t}(A)}\!\!\!\subHarmonicU(x)\dif x
\end{equation}
is a positive supermartingale. 
\end{theorem}
\begin{proof}
Denote $M_t:=\subHarmonicU(\proup{t})|D\proup{t}|$.
For every $0\le s\le t $ and $n\ge 1$ define 
\begin{equation}
S_n=\min_{i=1,\dots,d}\inf \left\{\tau\ge s:\left| \sum_{k=1}^d M_\tau \frac{\subHarmonicU_k(\proup{\tau})}{\subHarmonicU(\proup{\tau})}\sigma^{ki} \right|\ge n\right\}
\end{equation}
and integrate \eqref{M_supermartingale} on $[s,t\wedge S_n]$, 
\begin{equation}\label{equation_for_M}
M_{t\wedge S_n} = M_s + \int_s^{t\wedge S_n} M_\tau \frac{\nabla\subHarmonicU(\proup{\tau})}{\subHarmonicU(\proup{\tau})} \cdot \sigma\dif b
+ \int_s^{t\wedge S_n} M_\tau \frac{(\gsharpAdj \subHarmonicU)(\proup{\tau})}{\subHarmonicU(\proup{\tau})}\dif \tau.
\end{equation}
Since the stochastic integral has expectation zero and $\gsharpAdj\subHarmonicU\le 0$, then ${\bf E}[M_{t\wedge S_n}|\mathcal{F}_s]\le M_s$. Because $S_n\uparrow\infty$, an application of Fatou's lemma gives ${\bf E}[M_{t}|\mathcal{F}_s]\le M_s$.

For almost every Brownian path, $M_t$ is a positive and continuous map in the variable $x$, hence its integral is well defined. By the formula for the change of variables we have $\int_{\proup{t}(A)}\subHarmonicU=\int_A M_t$ and then using Fubini's theorem 
\begin{equation}
{\bf E}\Big[\int_{\proup{t}(A)}\!\!\!\subHarmonicU\ \Big|\mathcal{F}_s\Big]= {\bf E}\Big[\int_A M_{t}\ \Big|\mathcal{F}_s\Big]\le \int_A M_s = \int_{\proup{s}(A)}\!\!\!\subHarmonicU.
\end{equation}
\end{proof}
Doob's theorem guarantees that $\lim_{t\to\infty}\int_{\proup{t}(A)}\subHarmonicU$ exists almost surely. It is shown next that it is always possible to find many $u$'s such that the limit is zero. Condition \eqref{negativeDivergence}, as expected, will play a crucial role in the proof.

\begin{theorem}\label{imageOfAtendsToZero}
Let $\mathcal{W}=\{0<u\in C^{2}(\real{d}): u=\sqrt{z},\, z\in\mathcal{S}$\}.
Then for every $u\in \mathcal{W}$,
\begin{equation}
{\bf P}\Big[\lim_{t\to\infty}\int_{\proup{t}(A)}\!\!\!\subHarmonicU=0\Big]=1.
\end{equation}
\end{theorem}
\begin{proof}
By theorem \ref{princEigen_is_negative}, $\gsharpAdj\subHarmonicU\le -\frac{\lambda_c}{2}\subHarmonicU$.
Set $s=0$ in \eqref{equation_for_M} and take the expectation to obtain 
\begin{equation}
{\bf E}(M_{t\wedge S_n}) = M_0 
+ {\bf E}\!\!\int_0^{t\wedge S_n}\!\! M_\tau \frac{(\gsharpAdj \subHarmonicU)(\proup{\tau})}{\subHarmonicU(\proup{\tau})}\dif \tau\le M_0-\frac{\lambda_c}{2}\, {\bf E}\!\!\int_0^{t\wedge S_n}\!\! M_\tau\dif\tau.
\end{equation}
Sending $n\uparrow\infty$ and using Fatou's lemma we have 
\begin{equation}
{\bf E}(M_{t}) \le M_0 
-\frac{\lambda_c}{2}{\bf E}\int_0^{t} M_\tau\dif\tau.
\end{equation}

Next integrate over A,
\begin{equation}
{\bf E}\int_AM_{t}\le \int_A M_0-\frac{\lambda_c}{2}{\bf E} \int_0^t \Big(\int_A M_\tau\Big)\dif\tau,
\end{equation}
and observe that the second term on the RHS is monotone decreasing and bounded from below, hence ${\bf E}\int_0^\infty(\int_A M_\tau)
\dif\tau<\infty$. Therefore $\int_0^\infty(\int_A M_\tau)\dif\tau<\infty$ almost surely. Since 
$\int_{\proup{t}(A)}\subHarmonicU=\int_AM_t$ is integrable and convergent as $t\uparrow\infty$, the limit is necessarily zero.
\end{proof}

\section{Lyapunov exponent}
In this section the rate of convergence of $\int_{\proup{t}(A)}\subHarmonicU$ is analyzed to derive lower (upper) bound for the lim inf (lim sup) of 
\begin{equation}
\frac 1{t}\ln\!\!\!\int_{\proup{t}(A)}\!\!\!\subHarmonicU.
\end{equation}
The lower bound is easy to obtain.
\begin{theorem}\label{lyapunovLowerBound}
Let $0<u\in C^{2}(\real{d})$ and $M_t:=u(\proup{t})|D\proup{t}|$. If
\begin{equation}
\int_{\real{d}}\frac{\nabla u\cdot a\nabla u}{u^2\psi^2}<\infty\mbox{ and }
\int_{\real{d}}\frac{|\gsharpAdj u|}{u\psi^2}<\infty,
\end{equation}
then
\begin{equation}\label{lowerBoundIntegrand}
\lim_{t\to\infty}\frac 1{t}\ln{M_t}=-\frac 1{2}\int_{\real{d}}\frac{\nabla u\cdot a\nabla u}{u^2\psi^2}
+\int_{\real{d}}\frac{\gsharpAdj u}{u\psi^2}\quad\mbox{a.s.}
\end{equation}
and
\begin{equation}\label{lowerBound}
\liminf_{t\to\infty}\frac 1{t}\ln\!\!\!\int_{\proup{t}(A)}\!\!\!u\ge \mathrm{Leb}(A)\left[-\frac 1{2}\int_{\real{d}}\frac{\nabla u\cdot a\nabla u}{u^2\psi^2}
+\int_{\real{d}}\frac{\gsharpAdj u}{u\psi^2}\right]\quad\mbox{a.s.}
\end{equation}
\end{theorem}
\begin{proof} To begin with, recall from lemma \ref{lemma_uXdet(jac)} that
\begin{equation}
M_t=u(x)\mathrm{exp}\left(\int_0^t\frac{\nabla u(\proup{\tau})}{u(\proup{\tau})}\cdot\sigma\dif b -\frac{1}{2}\int_0^t \frac{(\nabla u\cdot a\nabla u)(\proup{\tau})}{u^2(\proup{\tau})}\dif \tau+\int_0^t\frac{(\gsharpAdj u)(\proup{\tau})}{u(\proup{\tau})}\dif \tau  \right).
\end{equation}
It is easy to show that the process
\begin{equation}
B=\int_0^t\frac{\nabla u(\proup{\tau})}{u(\proup{\tau})}\cdot\sigma\dif b
\end{equation}
is a standard 1-dimensional Brownian motion run with a random clock $\mathfrak{t}(t)=\int_0^t \frac{(\nabla u\cdot a\nabla u)(\proup{\tau})}{u^2(\proup{\tau})}\dif \tau$. By virtue of the ergodic theorem for recurrent diffusions \cite{tanaka} (see also \cite{pinsky}, Th. 4.9.5), 
$\lim_{t\to\infty}\frac{\mathfrak{t}(t)}{t}=\int_{\real{d}}\frac{(\nabla u\cdot a\nabla u)}{u^2\psi^2}$ almost surely, leading to
\begin{equation}
\lim_{t\to\infty}\frac{1}{t}\ln\left(M_t\mathrm{e}^{-B(\mathfrak{t})}\right)=-\frac 1{2}\int_{\real{d}}\frac{\nabla u\cdot a\nabla u}{u^2\psi^2}
+\int_{\real{d}}\frac{\gsharpAdj u}{u\psi^2}\quad\mbox{a.s.}
\end{equation}
To conclude the first half of the proof, it must be checked that almost surely
$\limsup_{t\to\infty}{\frac{B(\mathfrak{t})}{t}}=\liminf_{t\to\infty}{\frac{B(\mathfrak{t})}{t}}=0$. 
Since 
\begin{equation}
{\frac{B(\mathfrak{t})}{t}}={\frac{B(\mathfrak{t})}{\sqrt{2\mathfrak{t}\ln\ln{\mathfrak{t}}}}}\ \sqrt{\frac{2\ln\ln{\mathfrak{t}}}{\mathfrak{t}}}\ \frac{\mathfrak{t}}{t}
\end{equation}
and almost surely $\lim_{t\to\infty}\frac{\ln\ln{\mathfrak{t}}}{\mathfrak{t}}=0$, 
an application of the \v{H}in\v{c}in's law of the iterated logarithm proves it.

To prove \eqref{lowerBound} first apply the reversed Jensen's inequality for concave functions to obtain
$\ln\int_{A}M_t\ge\int_A\ln M_t$. Then Fatou's lemma permits to let $t\uparrow\infty$ under the integral sign, producing
\begin{equation}
\liminf_{t\to\infty}\frac 1{t}\ln\!\!\!\int_{\proup{t}(A)}\!\!\!u\ge\int_{A}\liminf_{t\to\infty}\frac 1{t}\ln M_t\quad\mbox{a.s.}
\end{equation}
Formula \eqref{lowerBoundIntegrand} supplies the rest of the proof.
\end{proof}

The following result identifies the upper bound. Condition \eqref{negativeDivergence} is critical.
\begin{theorem}\label{lyapunovUpperBound}
Let $0<\subHarmonicU\in C^{2}(\real{d})$ such that $\gsharpAdj\subHarmonicU\le -\frac{\lambda_c}{2}\subHarmonicU$, then
\begin{equation}
\limsup_{t\to\infty}\frac 1{t}\ln\!\!\!\int_{\proup{t}(A)}\!\!\!\subHarmonicU\le
-\frac{\lambda_c}{2}\quad\mbox{a.s.}
\end{equation}
\end{theorem}
\begin{proof}
First recall that theorem \ref{princEigen_is_negative} guarantees the existence of infinitely many smooth $\subHarmonicU>0$ such that $\gsharpAdj\subHarmonicU\le -\frac{\lambda_c}{2}\subHarmonicU$. Next define $\mathfrak{Z}_t:=e^{\lambda_ct/2}M_t$ and let $S_n$ be as in Section \ref{focusingflow}. It\^{o}'s Lemma gives
\begin{equation}
\mathfrak{Z}_{t\wedge S_n} \le \mathfrak{Z}_s + \int_s^{t\wedge S_n} \mathfrak{Z}_\tau \frac{\nabla\subHarmonicU(\proup{\tau})}{\subHarmonicU(\proup{\tau})} \cdot \sigma\dif b,
\end{equation}
permitting to conclude, after an application of Fatou's lemma, that $\mathfrak{Z}_t$ and $\int_A \mathfrak{Z}_t$ are both positive supermartingales. It follows that there must exist a finite random variable $V(b)$ such that $\sup_{t\ge 0}\int_A\mathfrak{Z}_t\le V(b)$, thus
\begin{equation}
\frac{1}{t}\ln\int_A M_t\le \frac{1}{t}\ln V(b)-\frac{\lambda_c}{2}\quad\mbox{a.s.}
\end{equation}
Taking the $\limsup_{t\to\infty}$ completes the proof.
\end{proof}
\section{An example: The Ornstein-Uhlenbeck process}
Here $\sigma$ is the identity matrix and $m(x)=-x$. To begin with, observe that the drift is generated by the potential
\begin{equation}
\psi(x)=e^{\frac{(x,x)^2}{2}}
\end{equation}
and since $\gupAdj (1/\psi^2)=0$, it follows that  
$1/\psi^2$ is the right invariant probability density (after, of course, a normalization).
Furthermore 
\begin{equation}\label{criticalCoefficient}
\nabla\cdot m=-d,
\end{equation}
so $\lambda_c=d$.
Now, simple computations show that
\begin{eqnarray}
 &  & \gsharpAdj \psi^2 =0,\\[3mm] 
 & & \gsharpAdj \psi = -\frac{1}{2\psi}\left(\nabla\psi\cdot\nabla\psi-\psi^2\nabla\cdot m\right),\label{gSharpAdjOnpsi}\\[3mm]
 & & \frac{\nabla\psi\cdot\nabla\psi}{\psi^4}=\frac{(x,x)}{\psi^2}.
 \end{eqnarray} 
Because the first term on the RHS of \ref{gSharpAdjOnpsi} is non-positive, using \ref{criticalCoefficient} we conclude $\gsharpAdj \psi\le -\frac{d}{2}\psi$.
Next choose $u=\psi$ and check that the conditions in theorems \ref{lyapunovLowerBound} and \ref{lyapunovUpperBound} are satisfied. Hence, 
\begin{equation}
\mathrm{Leb}(A)\left[-\frac{d}{2}-\int_{\real{d}}\frac{\nabla \psi\cdot a\nabla \psi}{\psi^4}\right]\le
\liminf_{t\to\infty}\frac 1{t}\ln\!\!\!\!\int_{\proup{t}(A)}\!\!\!\psi\le\limsup_{t\to\infty}\frac 1{t}\ln\!\!\!\!\int_{\proup{t}(A)}\!\!\!\psi\le
-\frac{d}{2}\quad\mbox{a.s.}
\end{equation}

\section{Conclusions}
This paper has studied the behavior of a stochastic flow generated by a positive recurrent diffusion in \real{d}. Even though the main theorems are proved under ``unattractive" conditions on the coefficients, the asymptotic focusing of $\int_{\proup{t}(A)}\subHarmonicU$ is likely guaranteed under weaker assumptions. There is evidence that only an integrable invariant measure is needed to make Theorems \ref{imageOfAtendsToZero}, \ref{lyapunovLowerBound}, \ref{lyapunovUpperBound} hold. The uniform Lipschitz assumption on the drift is purely technical: it guarantees that the explosion time $\mathfrak{e}(x,b)$ of \proup{t} verifies
\begin{equation}
{\bf P}[\mathfrak{e}(x,b)=\infty ,\mbox{ for all } x\in\real{d}]=1,
\end{equation}
which is slighty stronger than bare positive recurrence. Similarly,
the reversibility condition, critical to prove transience of \prosharp{t}, should not be necessary. Unfortunately, as R. Pinsky points out in \cite{pinsky}, the analysis 
of a non-symmetric generator in \real{d} is quite difficult.

The key property to exploit in future investigations seems to be the critical nature of the operator $\gsharpAdj$ which offers a broad set of superharmonic functions.  

The challenged reader may be interested in the following problem: Let $A$ be a compact set in $\real{d}$ and let $\proup{t}$ be a reversible recurrent diffusion in $\real{d}$ with smooth coefficients and integrable invariant density $\psi^{-2}$. Prove that  
\begin{equation}
\lim_{t\to\infty}\frac 1{t}\ln\!\!\!\int_{\proup{t}(A)}\!\!\!\psi^2= \mathrm{Leb}(A)\Big[-2\int_{\real{d}}\frac{\nabla \psi\cdot a\nabla \psi}{\psi^4}\Big]\quad\mbox{a.s.}
\end{equation}

\bibliography{mybibliography}
\end{document}